\documentclass[reqno,12pt]{amsart}

\usepackage{amsmath,amssymb}
\usepackage{amscd,amsfonts}
\usepackage{amsthm}
\usepackage{amsbsy}
\usepackage[ansinew]{inputenc}
\usepackage[T1]{fontenc}

\usepackage{hyperref}

\swapnumbers\theoremstyle{plain}

\newtheorem{thm}{Theorem}[section]
\newtheorem{prop}[thm]{Proposition}
\newtheorem{lem}[thm]{Lemma}
\newtheorem{cor}[thm]{Corollary}

\theoremstyle{definition}

\newtheorem{defi}[thm]{Definition}
\newtheorem{ex}[thm]{Example}

\newtheorem{rem}[thm]{Remark}

\DeclareMathOperator{\ef}{\mathbb{F}}
\DeclareMathOperator{\enn}{\mathbb{N}}

\DeclareMathOperator{\aaa}{\mathbb{A}}
\DeclareMathOperator{\af}{\mathbf{A}}
\DeclareMathOperator{\ddd}{\mathbb{D}}
\DeclareMathOperator{\df}{\mathbf{D}}

\DeclareMathOperator{\zed}{\mathbb{Z}}
\DeclareMathOperator{\Hom}{Hom}
\DeclareMathOperator{\Homf}{Hom_{\ef}}

\DeclareMathOperator{\Hod}{Hom_{\df}}

\DeclareMathOperator{\Ker}{Ker}

\DeclareMathOperator{\Mod}{Mod}
\DeclareMathOperator{\Modf}{Modf}

\DeclareMathOperator{\vect}{Vect}
\newcommand{\prf}{\noindent\textit{Proof. }}
\usepackage{color}

\begin{document}

\title{Discrete linear Algebraic Dynamical Systems}
\author{R. ANDRIAMIFIDISOA* and J. ANDRIANJANAHARY$^+$}
\begin{abstract}The vector space of the multi-indexed sequences over a field and the vector space of the sequences  with finite support are dual to each other, with respect to a \textit{scalar product}, which we used to define \textit{orthogonals} in these spaces. The closed subspaces in the first vector space are then the orthogonals of subsets in the second space. Using power series and polynomials, we prove that the \textit{polynomial operator in the shift} which U. Oberst and J. C. Willems  have introduced to define time invariant discrete linear dynamical systems is the functorial adjoint of the polynomial multiplication.  These results are generalized to the case of vectors of sequences and vectors of power series and polynomials. We end this paper by describing  discrete linear algebraic dynamical systems.\\

\end{abstract}
\maketitle
\noindent*{\tiny Corresponding author : Ramamonjy ANDRIAMIFIDISOA, Département de Mathématiques et Informatique, Faculté des Sciences, Université d'Antananarivo, BP 906, Madagascar, e-mail: ramamonjy.andriamifidisoa@univ-antananarivo.mg}\\
\noindent$^+${\tiny coauthor: Juanito ANDRIANJANAHARY, e-mail: juanitorabibisoa@yahoo.fr}
\section*{Introduction}Discrete algebraic dynamical systems theory essentially studies subsets $\boldsymbol{\mathcal{B}}$,  called \textit{behavior},  of the set of functions from a time set $\mathbb{T}$ (usually $\enn^r,\zed^r) $ to $\ef^l$, where $\ef$ is a field. According to  Oberst \cite{ob} and Willems \cite{wil1,wil2,wil3,wil5},  systems are linear, time invariant and  closed with respect to the topology of pointwise convergence. The purpose of this paper is to  find and prove these properties for the  mathematical structures  behind these systems,  and then to deduce these   properties for systems.
These results leads to an elegant presentation  and a characterization of  discrete linear algebraic dynamical systems. \\
This paper is organized as follows :

In Section \ref{section1}, we present the basic mathematical structures $\ef^{\omega}$ and $\ef^{(\omega)}$. Then we define scalar products and present the duality property. Finally, we  define orthogonalities in these spaces.\\
In Section \ref{section2} , we present the topology on $\aaa=\ef^{\enn^r}$. The main theorem of the section is Theorem \ref{o-basis}, which allows the construction  of an $0$-basis in $\aaa$.\\
In Section \ref{section3} we present Theorem \ref{closed}, which characterizes  the closed subspaces of $\aaa$.\\
In Section \ref{section4}, we generalize the preceding results to the sets $\aaa^l$ and $\ddd^l$. The main results of the section  are Theorem \ref{l-o-basis},  which allows the construction of an  $0$-basis in $\aaa^l$ and Theorem \ref{l-closed} which characterizes  the closed subspaces of $\aaa^l$.\\
In Section \ref{section5}, we introduce the vector spaces of polynomials $\df$, the vector spaces of formal powers series $\af$  and define orthogonalities in $\df^l$ and $\af^l$. The main result of the section is Theorem \ref{orth-closed} which gives the closed subsets of $\af^l$.   \\
In Section  \ref{section6},  we introduce the\textit{ shift operator} and the \textit{polynomial operator in the shift}. These lead to the \textit{polynomial-power series multiplication}, denoted by ``$\circ$'', the main tool for constructing discrete linear dynamical systems. The main results of the section are Theorem \ref{fundamental} and \ref{l-fundamental}, characterizing the polynomial and power series multiplication. \\
In the last section \ref{section7}, we define discrete algebraic dynamical systems, parting from our main result, Theorem \ref{autoregressive-orth-closed}, which presents various equivalent  properties of systems, and present one of their main properties, Theorem \ref{system-hom}.
\section{Duality of vector spaces}\label{section1}
\noindent Let $r\geqslant 1$ be an integer, $\leqslant$  the usual total ordering on the integers and $\leqslant_+$ the partial  ordering on $\enn^r$
defined by
 \begin{equation*}
   \alpha=(\alpha_1,\ldots,\alpha_r)\leqslant_+\beta=(\beta_1,\ldots,\beta_r)\Longleftrightarrow\alpha_i\leqslant\beta_i\quad\;\text{for}\; i=1,\ldots,r.
 \end{equation*}
Let $\psi_r$ be the mapping
\begin{align}\label{order-isom-1}
\begin{split}
  \psi_r:\enn &\longrightarrow \enn^r \\
  n&\longmapsto \psi_r(n)=\underbrace{(n,n,\ldots,n)}_{r\; \text{times}}.
\end{split}
\end{align}
We have\begin{equation}\label{order-isom-2}
    m<n \Longleftrightarrow \psi_r(m)<_+ \psi_r(n).
 \end{equation}
for $m,n\in\enn$. Let $\Omega(\enn)$ be the set
\begin{equation*}
\Omega(\enn)=\enn\bigcup_{r=1}^{+\infty}\{\enn^r\}=\{0,1,2,\ldots,n,\ldots\}\bigcup\{\enn, \enn^2,\enn^3,\ldots\}.
\end{equation*}
 The letter $\omega$ will denote an element of $\Omega(\enn)$ and  will be considered  as an \textit{ordinal}. A ordinal $\omega =n\in\enn$ is identified with the set $\{0,1,2,\ldots,n-1\}$. In this case, the notation $x\in\omega$ means that $x\in\{0,1,2,\ldots,n-1\}$. In the following, $\omega$ will denote a non-zero ordinal.  For $n\in\enn^*$, we denote by $\Delta_n$ the set
\begin{equation}\label{ord-inv}
 \Delta=\{\alpha\in\enn^r\;\vert\;\alpha\leqslant_+\psi_r(n)\}\subset\enn^r.
\end{equation}
Let $\ef$ be a field. All vector spaces  will be over $\ef$. For two vector spaces $E$ and $F$, we denote by $\Homf(E,F)$ the set of all linear mappings from $E$ to $F$ : it is again an $\ef$-vector space. For an ordinal $\omega$, we denote by $\omega$ the vector space of all mappings
\begin{align}\label{rcl}
\begin{split}
      y :\ef^{\omega}&\longrightarrow \ef \\
      \alpha                  & \longmapsto  y(\alpha)= y_\alpha.
    \end{split}
\end{align}
 Let $\ef^{(\omega)}$ be the vector subspaces of $x\in\ef^{\omega}$ with finite support :
\begin{equation*}
\ef^{(\omega)}=\{x\in\omega\;\vert\;\{\alpha\in\omega|\ x_\alpha\neq 0\}\quad\text{is a finite set}\}.
\end{equation*}
For $\alpha\in\omega$, let $\delta_\alpha$ be the element of $\ef^{(\omega)}$ defined by
\begin{equation}\label{delta}
\delta_\alpha(\beta)=\left\{
                       \begin{array}{ll}
                         1\quad\hbox{if}\quad\alpha= \beta,\\
                         0\quad\hbox{otherwise},
                       \end{array}
                     \right.
\end{equation}
for all $\beta\in\omega$. In other words, $\delta_\alpha(\beta)=\delta_{\alpha\beta}$ where $\delta_{\alpha\beta}$ is the Kronecker's symbol. Then $(\delta_\alpha)_{\alpha\in\omega}$ is an $\ef$-basis of $\ef^{(\omega)}$ and for all $x\in\ef^{(\omega)}$, we have
\begin{equation}\label{basis}
x=\sum_{\alpha\in \omega}x_\alpha\cdot \delta_\alpha.
\end{equation}
Now, an element $y\in\omega$ is defined by $(y_\alpha)_{\alpha\in\enn^r}$ where $y_\alpha=y(\alpha)$. We identify $y$ by the \textit{formal sum}
\begin{equation}
y=\sum_{\alpha\in\enn^r}y_\alpha\delta_\alpha,
\end{equation}
meaning that for all $\beta\in\enn^r$, the value of $y(\beta)$ is given by
\begin{equation}
y(\beta)=\sum_{\alpha\in\enn^r}y_\alpha\delta_{\alpha}(\beta).
\end{equation}
(This time, the sum is finite and given by $y_\beta$, thus well defined).\\
For $f\in\Hom_{\ef}(\ef^{(\omega)},\ef)$ and $x\in\ef^{(\omega)}$, using the equation \eqref{basis} above, we get
\begin{equation}\label{hom1}
 f(x)=\sum_{\alpha\in \omega}x_\alpha\cdot f(\delta_\alpha)
\end{equation}
Therefore $f$ is defined by  $(f(\delta_\alpha))_{\alpha\in\omega}\in\ef^{\omega}$. Conversely, the element  $f=(f_\alpha)_{\alpha\in\omega}$ of $\omega$ defines an element of $\Hom_{\ef}(\ef^{(\omega)},\ef)$ by  $\delta_\alpha\mapsto f_\alpha$ for all $\alpha\in\omega$, i.e
\begin{equation*}
   \Bigl( \forall x\in\omega\Bigr)\;\; \Bigr[ x\longmapsto \sum_{\alpha\in \omega}x_\alpha\cdot f_\alpha\Bigr].
\end{equation*}
The following proposition is fundamental:
\begin{prop}\label{scal-prod} The mapping
\begin{align}\label{def-scal}
    \begin{split}
    \langle -,- \rangle : \ef^{(\omega)}\times\ef^{\omega}&\longrightarrow \ef \\
      (x,f)                  & \longmapsto    \langle x,f \rangle=f(x)=\sum_{\alpha\in \omega}x_\alpha\cdot f(\delta_\alpha).
\end{split}
\end{align}
satisfies to the following properties:\\
(1) The homomorphism
\begin{align}\label{left}
\begin{split}
 \ef^{(\omega)}&\longrightarrow\Homf(\ef^{\omega},\ef)\\
 x&\longmapsto \left\{
                 \begin{array}{ll}
                   \langle x,-\rangle :&\omega \longrightarrow \ef \\
                   &f\longmapsto f(x)
                 \end{array}
               \right.
 \end{split}
\end{align}
and
\begin{align}\label{right}
\begin{split}
 \ef^{\omega}&\longrightarrow\Homf(\ef^{(\omega)},\ef)\\
 f&\longmapsto \left\{
                 \begin{array}{ll}
                   \langle -,f\rangle :&\ef^{(\omega)} \longrightarrow \ef \\
                   &x\longmapsto f(x)
                 \end{array}
               \right.
\end{split}
\end{align}
are injective .\\
(2) The  monomorphism \eqref{right} is an isomorphism of vector spaces.
\end{prop}
We say that $\langle -,- \rangle $ is a \textit{scalar product} and the vector spaces $\ef^{\omega}$ and $\ef^{(\omega)}$ are \textit{dual}.

For $P\subset\ef^{\omega}$ and $Q\subset\ef^{(\omega)}$, we define the orthogonal $P^{\perp}$ and $Q^{\perp}$ by
\begin{align}
  P^{\perp} &=\{x\in \ef^{(\omega)}\ |f(x)=0\ \forall f\in P\} \subset \ef^{(\omega)}\label{P-orth}, \\
  Q^{\perp} &=\{f\in \ef^\omega \ |f(x)=0\ \forall x\in Q\} \subset\ef^{\omega}\label{Q-orth}.
\end{align}
Note that $P^{\perp}$ (resp. $Q^{\perp}$ ) is a vector subspace of $\ef^{(\omega)}$ (resp. $\ef^{\omega}$). We end this section by stating the following lemma and its corollary, whose proofs are straightforward.
\begin{lem}\label{orth-corr}
Let  $P,P'$ be  subsets of $\ef^{\omega}$ and $Q,Q'$  subsets of $\ef^{(\omega)}$. Then
\begin{align}
  P\subset P'\Longrightarrow\ P^{\perp} \supset P'^{\perp}\quad\text{and}\quad Q\subset Q'\Longrightarrow Q^{\perp} \supset Q'^{\perp}\\
  P\subset P^{{\perp}{\perp}}\quad\text{and}\quad Q\subset Q^{{\perp}{\perp}}\label{2-orth}\\
  P^{{\perp}{\perp}{\perp}}=P^{\perp}\quad\text{and}\quad Q^{{\perp}{\perp}{\perp}}=Q^{\perp}.\label{3-orth}
\end{align}
\end{lem}
 \begin{cor}\label{gal-corr}
Consider the sets
 \begin{equation*}
    \mathcal{L}=\{P^{\perp}| P\subset\ef^{\omega} \} \text{ and } \mathcal{CL}=\{Q^{\perp}| Q\subset\ef^{(\omega)} \}.
 \end{equation*}
 Then the map
 \begin{align}\label{bij1}
 \begin{split}
      \mathcal{L}&\longrightarrow \mathcal{CL} \\
      P^{\perp}    & \longmapsto    P^{{\perp}{\perp}}
 \end{split}
 \end{align}
 is a bijection. Its inverse is
\begin{align}
 \begin{split}
\mathcal{CL}&\longrightarrow \mathcal{L} \\
      Q^{\perp}    & \longmapsto    Q^{{\perp}{\perp}}.
\end{split}
 \end{align}
 \end{cor}

\section{Topology on $\ef^{\enn^r}$}\label{section2}

By considering $\ef^{\omega}$ for $\omega=\enn^r$, we get the vector space $\aaa=\ef^{\enn^r}$. Using Proposition \ref{scal-prod}, $\aaa$ is in duality with $\ddd=\ef^{(\enn^r)}$  thanks to the scalar product \eqref{def-scal}. We therefore may consider the orthogonal $\perp$ defined in \eqref{P-orth} and \eqref{Q-orth}, on the subsets of $\aaa$ and  $\ddd$.\\

 Considering $\ef$  as a topological vector space with the\textit{ discrete topology} \cite{b2}, a fundamental system of neighborhoods of  $0$ ($0$-basis) is the set $\{\{0\}\}$. A sequence $(a_n)_{n\in\enn^r}$ of elements of $\ef$ converges  to an element $a\in\ef$ if there exists an integer $N\in\enn$ such that
 \begin{equation*}
   n\geqslant N\Longrightarrow a_n=a.
 \end{equation*}

The vector space  $\aaa$ is provided with the \textit{product topology} of that of the $\ef$'s and becomes a topological vector space too (\cite{b2,ob}). According to the definition of the product topology, we can present an $0$-basis of $\aaa$ which is a  modification of  what is mentioned in \cite{ob}. It uses the mapping $\psi_r$ defined in \eqref{order-isom-1}:
\begin{thm}\label{o-basis} An $0$-basis of  $\aaa$ is given by the family of sets $(V_n)_{n\in\enn^*}$ where
\begin{equation}\label{vn}
    V_n=\{y\in\aaa|\ y_\alpha=0 \text{ for } \alpha\leqslant_+\psi_r(n) \} \;\text{for}\;n\in\enn^*.
\end{equation}
It has  the property
\begin{equation*}
    V_1\supset V_2\supset V_3\supset\cdots \supset V_n\supset V_{n+1}\supset\cdots .
\end{equation*}
\end{thm}
\prf  We will use  some notations. If $I$ is a set of indices, $(E_i)_{i\in I}$ a family of sets indexed by $I$, and $A$ a set, then $\prod_{i\in I} A$ denotes   the cartesian product $\prod_{i\in I}E_i$ when $E_i=A$ for $i\in I$. It is also the set $A^I$ of the mappings from $I$ to $A$ :
\begin{equation*}
  A^I=\prod_{i\in I} A.
\end{equation*}
In particular, for $A=\ef$ and $I=\enn^r$, we have
\begin{equation*}
  \aaa=\ef^{\enn^r}=\prod_{\alpha\in\enn^r}\ef.
\end{equation*}
The sets $\mathcal{O}$ defined by
\begin{equation}\label{open}
  \mathcal{O}=\prod_{\alpha\in I}\{0\}\times\prod_{\alpha\in\enn^r\setminus I}\ef,
\end{equation}
  where the $I$'s are finite subset of $\enn^r$ are open sets of $\aaa$ containing $0$ (called \textit{elementary rectangles}). The other open subsets of $\aaa$ are the union of such subsets  (\cite{ob,b4}). We make the convention \(\mathcal{O}=\aaa\) if $I=\emptyset$.\\

If $n\in\enn^*$ and\ $y\in V_n$, we have $y_\alpha=0$ for $\alpha\leqslant_+\psi_r(n)$. Let $\mathcal{O}$ be the set defined by
\begin{equation*}
\mathcal{O}=\prod_{\alpha\leqslant_+\psi_r(n)}\{0\}\times\prod_{\alpha\nleqslant_+\psi_r(n)}\ef,
\end{equation*}
If $z\in\mathcal{O}$, it verifies $z_\alpha=0$ for $\alpha\leqslant_+\psi_r(n)$, so that $z\in V_n$.
It follows that
\begin{equation}\label{vois}
V_n\supset\mathcal{O}.
\end{equation}
Let $\Delta_n$ be the set defined by
\begin{equation}\label{delta-set}
 \Delta_n=\{\alpha\in\enn^r\;\vert\;\alpha\leqslant_+\psi_r(n)\}\subset\enn^r.
\end{equation}
We have
\begin{equation*}
 \mathcal{O}=\prod_{\alpha\in\Delta_n}\{0\}\times\prod_{\alpha\in\enn^r\setminus\Delta_n}\ef,
\end{equation*}
 so that  $\mathcal{O}$ is then of the form \eqref{open}, By \eqref{vois}, the set $V_n$ is  a neighborhood of $0$.\\

If $\mathcal{U}$ another neighborhood of $0$, there exists a finite subset $I$ of $\enn^r$ such that
\begin{equation*}
  \mathcal{O}=\prod_{\alpha\in I}\{0\}\times\prod_{\alpha\in\enn^r\setminus I}\ef\subset\mathcal{U}.
\end{equation*}

  Let $n$ be an integer which is strictly greater than the maximum of the coordinates of the $\alpha$'s in $I$:
\begin{equation*}
  \Big(\forall \alpha\in I\Bigr)\quad \Bigl[ \alpha<_+\psi_r(n)\Bigr].
\end{equation*}
   Using again the set $\Delta_n$ as in \eqref{delta-set}, we  have
\begin{equation*}
I\subset \Delta_n\quad\text{and}\;\enn^r\setminus\Delta_n\subset\enn^r\setminus I ,
\end{equation*}
so that
\begin{align*}
y\in V_n&\Longrightarrow y_\alpha=0\quad\text{for}\;\alpha\leqslant_+\psi_r(n),\\
&\Longrightarrow y\in\prod_{\alpha\leqslant_+\Delta_n}\{0\}\times\prod_{\alpha\in\enn^r\setminus \Delta_n} \ef\subset\prod_{\alpha\in I}\{0\}\times\prod_{\alpha\in\enn^r\setminus I}\ef=\mathcal{O}.
\end{align*}
It follows that $V_n\subset \mathcal{O}$ and therefore $V_n\subset\mathcal{U} $. We have shown that the family $(V_n)_{n\in\enn^*}$  consists of neighborhoods of $0$ such that any other neighborhood of $0$ contains an element of this family. It is obvious that the sequence $(V_n)_{n\in\enn^*}$is decreasing. \qed
\begin{cor}A sequence $(f_n)_{n\in\enn}$ of elements of $\aaa$ converges to an element $f\in\aaa$ if and only if for $\alpha\in\enn^r$, the sequence $(f_{n\alpha})_{n\in\enn}$  converges to $f_\alpha$ in $\ef$.
\end{cor}
\prf Suppose that the sequence $(f_n)_{n\in\enn}$ of elements of $\aaa$ converges to an $f\in\aaa$ and let $\alpha$ be an element of $\enn^r$. There  exists  $m\in\enn$ such that $\alpha\leqslant_+\psi_r(m)$. Let $n\in\enn^*$ such that $n\geqslant m$ and $V_n$ the corresponding  neighborhood of $0$ in $\aaa$. There exists an $N\in\enn$ such that
\begin{align}\label{f-conv}
\begin{split}
k\geqslant N&\Longrightarrow (f_k-f)\in V_n\\
  &\Longrightarrow f_{k\beta}=f_\beta\quad\text{for}\;\beta\leqslant_+\psi_r(n).
\end{split}
\end{align}
Applying  the second equation of \eqref{f-conv} for the case $\beta=\alpha$, we have
\begin{equation*}
k\geqslant N\Longrightarrow f_{k\alpha}= f_\alpha,
\end{equation*}
i.e.  the sequence $(f_{n\alpha})_{n\in\enn}$  converges to $f_\alpha$ in $\ef$.\\

Conversely, suppose that  for $\alpha\in\enn^r$,  the sequence $(f_{n\alpha})_{n\in\enn}$  converges to $f_\alpha$ in $\ef$. Given $\alpha\in\enn$, there exists an $N\in\enn$ such that
 \begin{align}\label{a-conv}
k\geqslant N\Longrightarrow f_{k\alpha}-f_\alpha=0.
\end{align}
Let $V_n$ be a neighborhood of $0$ in $\aaa$ and $\Delta_n$ the set
\begin{equation}
 \Delta_n=\{\alpha\in\enn^r\;\vert\;\alpha\leqslant
 _+\psi_r(n)\}\subset\enn^r.
\end{equation}

Let $m$ be the cardinality of $\Delta_n$ and $\alpha_1, \ldots,\alpha_m$ the elements of $\Delta_n$. Applying \eqref{a-conv} for each of the elements of $\Delta_n$, there exist $N_1,\dots, N_m\in\enn$ such that
\begin{align*}
k\geqslant N_i\Longrightarrow f_{k{\alpha_i}}-f_{\alpha_i}=0.
\end{align*}
for $i=1\ldots,m$. Taking  $N=\max\{N_1,\ldots,N_m\}$,   we  have
\begin{align*}
k\geqslant N\Longrightarrow f_{k{\alpha_i}}-f_{\alpha_i}=0
\end{align*}
for $i=1\ldots,m$. In other words,
\begin{align*}
k\geqslant N&\Longrightarrow f_{k\alpha}=f_\alpha\quad\text{for}\;\alpha\leqslant_+\psi_r(n).
\end{align*}
We then have shown that
\begin{align*}
k\geqslant N&\Longrightarrow (f_k-f)\in V_n.
\end{align*}
Thus $(f_n)_{n\in\enn}$ converges to $f$ in $\aaa$. \qed\\

The topology of $\aaa$ is then  the topology of the \textit{pointwise convergence}.

\section{Closed subspaces of the vector space of multi-indexed sequences over a field}\label{section3}
We are now going to investigate the closed subspaces of $\aaa$.
\begin{thm}\label{closed}
An $\ef$-vector subspace  $V$  of $\aaa$ is closed if and only if
there exists $G\in\ddd$ such that
\begin{equation}\label{fsevclos}
    V=G^{\perp}.
\end{equation}
\end{thm}
\prf Let $G$ be a non-empty subset of $\ddd$. We have to show that $G^{\perp}$ is closed in $\aaa$. Let $(f_n)_{n\in\enn}$ a sequence in $G^{\perp}$ which converges to $f\in\aaa$, with respect to the topology of $\aaa$. Given $x\in G$, the following property holds:  for  $k\in\enn$, there exists $N\in\enn$ such that for $n\geqslant N$, one has $(f_n-f)\in V_k$, i.e. $(f_n-f)_\alpha=0$ whenever $\alpha\leqslant_+\psi_r(k)$. Since
\begin{align}
\begin{split}
    f(x)&=\sum_{\alpha\leqslant_+ \psi_r(k)}x_\alpha f_\alpha + \sum_{\alpha\nleqslant_+\ \psi_r(k)}x_\alpha f_\alpha\\
    &=\sum_{\alpha\leqslant_+ \psi_r(k)}x_\alpha (f_n)_\alpha + \sum_{\alpha{\nleqslant_+} \psi_r(k)}x_\alpha f_\alpha.
\end{split}
\end{align}
Since $f_n$ is with finite support, we may choose $k$ sufficiently large so that
\begin{equation*}
 0= f_n(x)=\sum_{\alpha\leqslant_+\psi_r(k)}x_\alpha f_{n\alpha}.
\end{equation*}
 It follows that\begin{equation*}
    f(x)=\sum_{\alpha\nleqslant_+ \psi_r(k)}x_\alpha f_\alpha.
\end{equation*}
But $x$ being with finite support, we may increase $k$, if necessary,  so that $x_\alpha=0$ for  $\alpha$ verifying $\alpha\nleqslant_+\psi_r(k)$, which implies $f(x)=0$. Since  it is true for an arbitrary $x\in G$, we finally have $f\in G^{\perp}$.\\
Conversely, suppose that $V$ is closed in $\aaa$. We will show that
\begin{equation}\label{double-orth}
V^{{\perp}{\perp}}=V=\overline{V},
\end{equation}
where $\overline{V}$ is the closure of $V$ with respect to the topology of $\aaa$. It suffices to show the non-trivial inclusion $V^{{\perp}{\perp}}\subset V $ (see \eqref{2-orth}). Let $\{f_{\lambda}\;\vert\;\lambda\in \Lambda\}$ a  generating set  of $V$ :
\begin{align*}
\begin{split}
              f_{\lambda}:\ \enn^r &\longrightarrow \ef \\
              \alpha & \longmapsto  f_{\lambda\alpha}.
\end{split}
\end{align*}
Let $q=(q_\alpha)_{\alpha\in\enn^r}\in V^{{\perp}{\perp}}$, $n\in \enn$ and $G_n$ be the  finite-dimensional vector subspace of $\ddd$  defined by
\begin{equation}\label{gn}
G_n=\oplus_{\alpha\in\Delta_n}\ef \delta_\alpha
\end{equation}
where $\Delta_n=\{\alpha\in\enn^r\;\vert\;\alpha\leqslant_+\psi_r(n)\}$ ($G_n$ is the subspace of $\ddd$ generated by $\{\delta\alpha |\ \alpha\leqslant_+\psi_r(n)\}$). From classical linear algebra, we have the isomorphism
\begin{align*}
\begin{split}
  \ef^{\Delta_n}& \longleftrightarrow \Hom_{\ef}(G_n,\ef) \\
  x=(x_\alpha)_{\alpha\in\Delta_n}                      & \longleftrightarrow \left\{
                                                                       \begin{array}{ll}
                                                                       \varphi:&\ G_n  \rightarrow \ef \\
                                                                      & \delta_\alpha\mapsto x_\alpha\ .
                                                                       \end{array}
                                                        \right.
\end{split}
\end{align*}
For $\beta\in\Delta_n$, we define the element $\gamma_\beta\in \Hom_{\ef}(G_n, \ef)$ by
\begin{align*}
    \begin{split}
       \gamma_\beta:\ G_n & \longrightarrow \ef \\
       \delta\alpha     & \longmapsto    \gamma_\beta(\delta\alpha)=\delta_{\beta\alpha},
     \end{split}
\end{align*}
where $\delta_{\beta\alpha}$ is the Kronecker's symbol. The family  $\{\gamma_\beta\ |\ \beta\in  \Delta_n\}$ is an $\ef$-basis of $\Hom_{\ef}(G_n,\ef)$. Therefore, if $\varphi\in \Hom_{\ef}(G_n, \ef)$ then $\varphi$ may be written as
\begin{equation}\label{phi}
\varphi=
\sum_{\beta\in \Delta_n} \varphi_\beta \gamma_\beta
\end{equation}
with $\varphi_\beta\in \ef $. \\
Again, we have the isomorphism :
\begin{align}\label{isom2}
\begin{split}
\Phi: \ef^{\Delta_n}&\longleftrightarrow \Hom_{\ef}(\Hom_{\ef}(G_n,\ef),\ef)\\
  y=(y_\alpha)_{\alpha_\in\Delta_n}&\longmapsto\left\{
                                                        \begin{array}{ll}
                                                         \Phi(y):\Hom_{\ef}(G_n,&\ef)\rightarrow \ef \\
                                                           &\gamma_\alpha \mapsto \Phi(y)(\gamma_\alpha)=y_\alpha.
                                                        \end{array}
                                                      \right.
\end{split}
\end{align}
Consider the restrictions
\begin{align}\label{lambda-n}
\begin{split}
f_\lambda|_{G_n}&=f_\lambda^{(n)}\quad\text{with}\;f_{\lambda}^{(n)}(\delta_\alpha)=f_{\lambda\alpha},\\
q|_{G_n}&=q_n\quad\text{with}\;q(\delta_\alpha)=q_\alpha
\end{split}
\end{align}
for $\alpha\in\Delta_n$. We will need the following lemma :
\begin{lem}\label{lemme}
$q_n\in \langle (f_\lambda^{(n)})_{\lambda\in \Lambda}\rangle$.
\end{lem}
\prf We have $q_n\in \Hom_{\ef}(G_n,\ef)$; suppose that $q_n\notin\langle (f_\lambda^{(n)})_{\lambda\in \Lambda}\rangle$. In thi case, we know that there exists   $\Theta\in  \Hom_{\ef}(\Hom_{\ef}(G_n,\ef),\ef)$ with $\Theta(f_\lambda^{(n)})=0$ for $\lambda\in \Lambda$ and $\Theta(q_n)=1$. By the isomorphism  \eqref {isom2}, there exists $y=(y_\alpha)_{\alpha_\in  \Delta_n}\in \ef^{\Delta_n}$ such that $\Phi(y)=\Theta$, i.e
$\Theta(\gamma_\alpha)=y_\alpha$ for  $\alpha\in\Delta_n\ $. Take
\begin{equation*}
g=\sum_{\alpha\in \Delta_n}y_\alpha \delta_\alpha \ \in G_n.
\end{equation*}
 Then
\begin{align}
\begin{split}
    f_\lambda^n(g)= f_{\lambda|Gn}(g)&=f_\lambda(\sum_{\alpha\in  \Delta_n}y_\alpha \delta_\alpha)\\&=\sum_{\alpha\in  \Delta_n}f_\lambda (y_\alpha \delta_\alpha)=\sum_{\alpha\in  \Delta_n}y_\alpha f_\lambda (\delta_\alpha)=\sum_{\alpha\in  \Delta_n}y_\alpha f_{\lambda\alpha}\\
    &=\sum_{\alpha\in  \Delta_n}\Theta(\gamma_\alpha)f_{\lambda\alpha}\\
    &=\sum_{\alpha\in \Delta_n}\Theta(f_{\lambda\alpha}\gamma_\alpha)=\Theta(\sum_{\alpha\in  \Delta_n}f_{\lambda\alpha}\gamma_\alpha)=\Theta(f_\lambda^{(n)})=0.
\end{split}
\end{align}

We then get
\begin{equation*}
(\forall \lambda\in\Lambda) \quad \bigg[f_\lambda (g)=0\bigg],
\end{equation*}
which proves that $g\in V^{\perp}$. But we also have
\begin{align*}
q(g)&=q|_{G_n}(g)=q_n(g)=q_n(\sum_{\alpha\in \Delta_n}y_\alpha \delta_\alpha)\\
&=\sum_{\alpha\in\Delta_n}y_\alpha q_{\alpha} = \sum_{\alpha\in\Delta_n}\Theta (\gamma_\alpha) q_{\alpha}=\Theta(\sum_{\alpha\in\Delta_n}\gamma_\alpha q_{\alpha}),
\end{align*}
and by \eqref{phi},
\begin{equation*}
  q_n=\sum_{\alpha\in\Delta_n}\gamma_\alpha q_{n_{\alpha}}=\sum_{\alpha\in\Delta_n}\gamma_\alpha q_\alpha,
\end{equation*}
so that
\begin{equation*}
q(g) =\Theta(q_n)=1.
\end{equation*}
 This implies that $g\notin V^{\perp\perp\perp}=V^\perp$, which is a contradiction. We conclude that necessarily   $q_n\in \langle (f_\lambda^{(n)})_{\lambda\in \Lambda}\rangle$.\qed\\

\noindent\textit{Proof of theorem \ref{closed}} (continued). For  $n \in \enn^*$, there is then a family $(\mu _{\lambda}^{(n)})_{\lambda\in \Lambda}$ with finite support such that
 \begin{equation*}
    q|_{G_n}=q_n=\sum_{\lambda\in \Lambda} \mu _{\lambda}^{(n)} f_\lambda^{(n)}.
\end{equation*}
Consider  the element $q'_n=\sum_{\lambda\in \Lambda} \mu _{\lambda}^{(n)} f_\lambda\in V$, which is an extension of $q_n$ on $\aaa$. For $\alpha \in \Delta_n$, we have, by \eqref{lambda-n}
 \begin{equation*}
  q(\alpha)=q'_n(\alpha),
\end{equation*}
so that $q-q'_n\in V_n$. Therefore,  the sequence  $(q'_n)_{n\in \enn^*}$ converges to $q$ in $\aaa$. Finally, we have  $q\in \overline{V}$. Hence $V^{{\perp}{\perp}}\subset V$ and the equality holds. Taking $G=V^{\perp}$, we get $G^{\perp}=V^{{\perp}{\perp}}=V$ and the theorem is proved.\qed

\section{Generalization to the  case of vectors of multi-indexed sequences}\label{section4}
We use the notations in the preceding sections. For an integer $l\geqslant 1$, we denote by $\aaa^l$ the set of column-vectors of elements of $\aaa$ with $l$ rows and $\ddd^l$  the set of row-vectors   of elements of $\ddd$ with $l$ columns :
\begin{align*}
&\aaa^l =\Bigg\{ w=\begin{pmatrix} w_1 \\
\vdots \\
w_l \\
\end{pmatrix} \;\Bigg|\; w_i\in\aaa\quad \text{for}\;i=1,\ldots,l  \Bigg\},\\&\ddd^l =\Bigg\{ d=\begin{pmatrix}d_1 ,&\ldots ,& d_l\;
\end{pmatrix}\;\Bigg|\;d_i\in\ddd \;\text{for}\;i=1,\ldots,l  \Bigg\}.
\end{align*}
For $w\in\aaa^l$, with the components $w_i\in\aaa$, we write
\begin{equation*}\label{d}
  w_i=(w_{i\alpha})_{\alpha\in\enn^r}
\end{equation*}
and for $d\in\ddd^l$, with the components $d_i\in\ddd$, we write
\begin{equation*}\label{w}
  d_i=(d_{i\alpha})_{\alpha\in\enn^r}
\end{equation*}
with $d_{i\alpha}=0$ except for a finite number of $\alpha$'s. As in Proposition \ref{scal-prod}, we define the following scalar product :
\begin{prop}\label{l-scalar-prod}The $\ef$-bilinear mapping
  \begin{align}
  \begin{split}
  \langle-,-\rangle:\ddd^l\times\aaa^l&\longrightarrow\ef,\\
  (d,w)&\longmapsto \langle d,w\rangle = \sum_{i=1}^{l}(\sum_{\alpha\in\enn^{r}}d_{i\alpha}\cdot w_{i\alpha})
  \end{split}
  \end{align}
  is a scalar product.
  \end{prop}
We then have
the vector spaces isomorphism
\begin{align}\label{l-isom}
\aaa^l&\cong\Homf(\ddd^l,\ef)\\
w&\longmapsto\langle -,w\rangle.
\end{align}

 As in \eqref{P-orth} and \eqref{Q-orth}, we  define \textit{orthogonals}:
\begin{defi}For $P\subset\aaa^l$ and $Q\subset\ddd^l$,  the orthogonals $P^{\perp}$ and $Q^{\perp}$ are
\begin{align}
  P^{\perp} &=\{d\in \ddd^l\; |\;\langle d,w\rangle=0\ \forall
  w\in P\} \subset \ddd^l\label{l-P-orth}, \\
  Q^{\perp} &=\{w\in \aaa^l\; |\;\langle d,w\rangle=0\ \forall d\in Q\} \subset\aaa^l\label{l-Q-orth}.
\end{align}
\end{defi}
For $P\subset\aaa^l $, the set $P^{\perp}$ is a vector subspace of $\ddd^l$ and for $Q\subset\ddd^l$, the set $Q^\perp$ is a vector subspace of $\ddd^l$. We also have the following results, as in Lemma \ref{orth-corr} and Corollary \ref{gal-corr}:\\

\begin{lem}\label{l-orth-corr}
Let , $P,P'$ be  subsets of $\aaa^l$ and $Q,Q'$  subsets of $\ddd^l$. Then
\begin{align}
  P\subset P'\Longrightarrow\ P^{\perp} \supset P'^{\perp}\quad\text{and}\quad Q\subset Q'\Longrightarrow Q^{\perp} \supset Q'^{\perp}\label{1-orth}\\
  P\subset P^{{\perp}{\perp}}\quad\text{and}\quad Q\subset Q^{{\perp}{\perp}}\label{2-orth}\\
  P^{{\perp}{\perp}{\perp}}=P^{\perp}\quad\text{and}\quad Q^{{\perp}{\perp}{\perp}}=Q^{\perp}.\label{3-orth}
\end{align}
\end{lem}

\begin{cor}
Consider the sets
 \begin{equation*}
    \mathcal{L}=\{P^{\perp}| P\subset\aaa^l \} \text{ and } \mathcal{CL}=\{Q^{\perp}| Q\subset\ddd^l\}.
 \end{equation*}
 Then the map
 \begin{align}\label{bij1}
 \begin{split}
      \mathcal{L}&\longrightarrow \mathcal{CL} \\
      P^{\perp}    & \longmapsto    P^{{\perp}{\perp}}
 \end{split}
 \end{align}
 is a bijection. Its inverse is
\begin{align}
 \begin{split}
\mathcal{CL}&\longrightarrow \mathcal{L} \\
      Q^{\perp}    & \longmapsto    Q^{{\perp}{\perp}}.
\end{split}
 \end{align}
 \end{cor}

As in section \ref{section2}, the field $\ef$  is provided with the discrete topology. The topology on $\ef^l$ is the product topology, therefore is also the discrete topology. A sequence $(a_n)_{n\in\enn}$ of elements of $\ef^l$ converges to an element $a\in\ef^l$ if there exists an $N\in\enn$ such that
\begin{equation}\label{l-F-conv}
  n\geqslant N\Longrightarrow a_n=a.
\end{equation}

 The topology on  $\aaa$ is the product of the discrete topology on $\ef$, with the $0$-basis given by Theorem \ref{o-basis} and the topology on $\aaa^l$ is the product of the topology on $\aaa$.
\begin{thm}\label{l-o-basis} An $0$-basis of  $\aaa^l$ is given by the family of sets $(\mathcal{O}_N)_{N\in\enn^*}$ where
\begin{align}\label{vn}
\begin{split}
    \mathcal{O}_N=\Bigg\{w=\begin{pmatrix}w_1 \\
\vdots\\
w_l \\
\end{pmatrix}\in\aaa^l\Bigg|\ w_{\alpha}&=0 \text{ for $\alpha\leqslant_+\psi_r(N)$ } \Bigg\} \\
&\text{for}\;N\in\enn^*.
\end{split}
\end{align}
It has  the property
\begin{equation*}
    \mathcal{O}_1\supset \mathcal{O}_2\supset \mathcal{O}_3\supset\cdots \supset \mathcal{O}_n\supset \mathcal{O}_{n+1}\supset\cdots .
\end{equation*}
\end{thm}
\prf \ According to the definition of the product topology on $\aaa^l$, the sets $V$ defined by
\begin{equation}\label{rect}
  V=\prod_{i=1}^lV_{i}
\end{equation}
  where the \(V_i\) are  elementary rectangles  of $\aaa$ are open sets of $\aaa$ containing $0$ and the other open sets of $\aaa^l$ containing $0$ are the union of such subsets  (\cite{b2,ob,b4}). \\
Using Theorem \ref{o-basis}, the set \(V_i\) is of the form
\begin{equation}\label{rect-comp}
  V_i=\prod_{\alpha\leqslant_+ \psi_r(n_i)}\{0\}\times\prod_{\alpha\nleqslant_+\psi_r( n_i)}\ef,
\end{equation}
for $i=1,\ldots,l$, where $n_i\in\enn^*$.\\

Using these notations, we have
\begin{equation*}
  \mathcal{O}_N=\prod_{i=1}^lV_N=V_N^l,
\end{equation*}
 so that the sets  $\mathcal{O}_N$ are indeed open sets of $\aaa^l$ containing $0$ (they are even elementary rectangles).\\
Now, let $V$ an elementary rectangle of $\aaa^l$ containing $0$, as in \eqref{rect}, where $V_i$ is given by \ref{rect-comp}. Let $N$ be the maximum of the $n_i$ :
\begin{equation*}
  N=\max\{n_i\;|\;i=1,\ldots,l\}.
\end{equation*}
Then $V_N^l\subset V$. Indeed, if $w\in V_N^l$, then  $w_{i\alpha}=0$ for $\alpha\leqslant_+\psi_r(N)$. But
\begin{equation*}
  \alpha\leqslant_+\psi_r(n_i)\Longrightarrow\alpha\leqslant_+\psi_r(N),
\end{equation*}
 so that $w_{i\alpha}=0$  whenever $\alpha\leqslant_+\psi_r(n_i)$. It follows that $w_i\in V_i$ for $î=1,\ldots,l$, i.e. $w\in V$.\\
We then have shown that the family of sets \(\mathcal{O}_N\) are open sets containing $0$, thus neighborhoods of $0$ in $\aaa^l$ such that every other neighborhood of $0$ contains an element of this family. In other terms, the family \(\mathcal{O}_{N\in\enn^*}\) is an $0$-base of $\aaa^l$.\qed

Let $(w_n)_{n\in\enn}$ a sequence of elements of $\aaa^l$ which converges to an element $w\in\aaa^l$, with $w_n=(w_{n\alpha})_{\alpha\in\enn^r}$ where $w_{n\alpha}\in\ef^l$ and $w=(w_\alpha)_{\alpha\in\enn^r}$, where $w_\alpha\in\ef^l$. Fix  $\alpha\in\enn^r$; there exists $k\in\enn$ such that $\alpha\leqslant_+\psi_r(k)$.   Fix $N\in\enn$ with $N\geqslant k$; there exists $M\in\enn$ such that
\begin{align*}
  n\geqslant M&\Longrightarrow w_n-w\in \mathcal{O}_N\\
    &\Longrightarrow w_{n\beta}=w_\beta\quad\text{for}\; \beta\leqslant_+\psi_r(N)\\
    &\Longrightarrow w_{n\alpha}=w_\alpha\quad\text{since}\;\;\alpha\leqslant_+\psi_r(k)\leqslant\psi_r(N).
\end{align*}
It follows that $(w_{n\alpha})_{n\in\enn}$ converges to $w_\alpha$ in $\ef^l$.\\
Conversely, suppose that for $\alpha\in\enn^r$, the sequence $(w_{n\alpha})_{n\in\enn}$ of $\ef^l$ converges to $w_\alpha$ in $\ef^l$.  Let $V_N$ be an element of the $0$-basis of $\aaa^l$, as in \ref{o-basis}. Given $\alpha\in\enn$, there exists an $M\in\enn$ such that
 \begin{align}\label{a-l-conv}
k\geqslant M\Longrightarrow w_{k\alpha}-w_\alpha=0.
\end{align}
As in \ref{delta-set}, let $\Delta_N$ be the set
\begin{equation}
 \Delta_N=\{\alpha\in\enn^r\;\vert\;\alpha\leqslant
 _+\psi_r(N)\}\subset\enn^r.
\end{equation}

Let $m$ be the cardinality of $\Delta_N$ and $\alpha_1, \ldots,\alpha_m$ the elements of $\Delta_N$. Applying \eqref{a-l-conv} for each of the elements of $\Delta_N$, there exist $M_1,\dots, M_m\in\enn$ such that
\begin{align*}
k\geqslant M_i\Longrightarrow w_{k{\alpha_i}}-w_{\alpha_i}=0.
\end{align*}
for $i=1\ldots,m$. Taking  $N=\max\{M_1,\ldots,M_m\}$,   we  have
\begin{align*}
k\geqslant N\Longrightarrow w_{k{\alpha_i}}-w_{\alpha_i}=0
\end{align*}
for $i=1\ldots,m$. In other words,
\begin{align*}
k\geqslant N&\Longrightarrow w_{k\alpha}=w_\alpha\quad\text{for}\;\alpha\leqslant_+\psi_r(N).
\end{align*}
We then have shown that
\begin{align*}
k\geqslant N&\Longrightarrow w_k-w\in V_N.
\end{align*}
Thus $(w_n)_{n\in\enn}$ converges to $w$ in $\aaa^l$. \qed\\

The topology of $\aaa^l$ is then also  the topology of the \textit{pointwise convergence}.\\
We are now going to investigate the closed subspaces of $\aaa^l$.
\begin{thm}\label{l-closed}
An $\ef$-vector subspace  $V$  of $\aaa^l$ is closed if and only if
there exists $G\in\ddd^l$ such that

\begin{equation}\label{fsevclos}
    V=G^{\perp}.
\end{equation}
\end{thm}

\prf Let $G$ be a non-empty subset of $\ddd^l$. We have to show that $G^{\perp}$ is closed in $\aaa^l$. Let $(w_n)_{n\in\enn}$ a sequence in $G^{\perp}$ which converges to $w\in\aaa^l$, with respect to the topology of $\aaa^l$. Write $w_n=(w_{ni})_{1\leqslant i\leqslant l}\in\aaa^l$ for $i\in\enn$ and $w=(w_i)_{1\leqslant i\leqslant l}\in\aaa^l$ Given $d\in G$, we have $\langle d,w_n\rangle=0$ for $n\in\enn$. For   $k\in\enn$, there exists $N\in\enn$ such that for $n\geqslant N$, one has $(w_n-w)\in V_k$, i.e. $(w_n-w)_\alpha=0$ whenever $\alpha\leqslant_+\psi_r(k)$. According to the notations in \ref{l-scalar-prod}, we have
\begin{align}
\begin{split}
  \langle d,w\rangle&=\sum_{i=1}^l\Big(\sum_{\alpha\leqslant_+ \psi_r(k)}d_{i\alpha} w_{i\alpha} + \sum_{\alpha\nleqslant_+\ \psi_r(k)}d_{i\alpha} w_{i\alpha}\Bigr)\\
    &=\sum_{i=1}^l\Bigl(\sum_{\alpha\leqslant_+ \psi_r(k)}d_{i\alpha} (w_{ni})_\alpha\Bigr) + \sum_{i=1}^l\Bigl(\sum_{\alpha{\nleqslant_+} \psi_r(k)}d_{i\alpha} w_{i\alpha}\Bigr).
\end{split}
\end{align}
Since $w_n$ is with finite support, we may choose $k$ sufficiently large so that
\begin{equation*}
 0= \langle d,w_n\rangle=\sum_{i=1}^l\Bigl(\sum_{\alpha\leqslant_+\psi_r(k)}d_{i\alpha}w_{ni\alpha}\Bigr).
\end{equation*}
 It follows that\begin{equation*}
    \langle d,w\rangle=\sum_{i=1}^l\Bigl(\sum_{\alpha\nleqslant_+ \psi_r(k)}d_{i\alpha} w_{i\alpha}\Bigr).
\end{equation*}
But $w$ being with finite support, we may increase $k$, if necessary,  so that $w_\alpha=0$ for  $\alpha$ verifying $\alpha\nleqslant_+\psi_r(k)$, which implies $\langle d,w\rangle=0$. Since  it is true for an arbitrary $d\in G$, we finally have $w\in G^{\perp}$.\\
Conversely, suppose that $V$ is closed in $\aaa^l$. We will show that
\begin{equation}\label{double-orth}
V^{{\perp}{\perp}}=V=\overline{V},
\end{equation}
where $\overline{V}$ is the closure of $V$ with respect to the topology of $\aaa^l$. It suffices to show the non-trivial inclusion $V^{{\perp}{\perp}}\subset V $ (see \eqref{2-orth}). Let $\{w_{\lambda}\;\vert\;\lambda\in \Lambda\}$ a  generating set  of $V$ :
\begin{align*}
\begin{split}
              w_{\lambda}:\ \enn^r &\longrightarrow \ef \\
              \alpha & \longmapsto  w_{\lambda\alpha}.
\end{split}
\end{align*}
Let $q=(q_\alpha)_{\alpha\in\enn^r}\in V^{{\perp}{\perp}}$, $n\in \enn$ and $G_n$ be the  finite-dimensional vector subspace of $\ddd^l$  defined by
\begin{equation}\label{gn}
G_n=\oplus_{(\alpha,i)\in\Delta_n}\ef \delta_\alpha e_i
\end{equation}
where $\Delta_n=\{\alpha\in\enn^r\;\vert\;\alpha\leqslant_+\psi_r(n)\}\times\{1,\ldots,l\}$, $\delta_\alpha=(\delta_{\alpha\beta})_{\beta\in\enn^r}$,  $\delta_{\alpha\beta}$  the Kronecker' symbol and
\begin{equation*}
  e_i=\underbrace{(0,\ldots,1,\ldots,0)}_{1\;\text{at the\;}k\;\text{-th \;position}},\quad \delta_\alpha e_i=\underbrace{(0,\ldots,\delta_\alpha,\ldots,0)}_{\delta_\alpha \;\text{at the}\; k\;\text{-th position}}\in\ddd^l
 \end{equation*}

 ($G_n$ is the subspace of $\ddd^l$ generated by $\{\delta_\alpha e_{i}|\ (\alpha,i)\in\Delta_n\}$). From classical linear algebra, we have the isomorphism
\begin{align*}
\begin{split}
  \ef^{\Delta_n}& \longleftrightarrow \Hom_{\ef}(G_n,\ef) \\
  x=(x_{(\alpha,i)})_{(\alpha,i)\in\Delta_n}                      & \longleftrightarrow \left\{
                                                                       \begin{array}{ll}
                                                                       \varphi:&\ G_n  \rightarrow \ef \\
                                                                      &\delta_ \alpha e_i\mapsto x_{(\alpha,i)}\ .
                                                                       \end{array}
                                                        \right.
\end{split}
\end{align*}
For $(\beta,j)\in\Delta_n$, we define the element $\gamma_{(\beta,i)}\in \Hom_{\ef}(G_n, \ef)$ by
\begin{align*}
    \begin{split}
       \gamma_{(\beta,i)}:\ G_n & \longrightarrow \ef \\
       \delta_\alpha e_i     & \longmapsto    \gamma_{(\beta,j)}(\delta_\alpha e_i)=\delta_{(\beta,j)(\alpha,i)},
     \end{split}
\end{align*}
where $\delta_{(\beta,j)(\alpha,i)}$ is the Kronecker's symbol. The family  $\{\gamma_{(\beta,j)}\ |\ (\beta,j)\in  \Delta\}$ is an $\ef$-basis of $\Hom_{\ef}(G_n,\ef)$. Therefore, if $\varphi\in \Hom_{\ef}(G_n, \ef)$ then $\varphi$ may be written as
\begin{equation}\label{phi}
\varphi=
\sum_{{(\beta,j)}\in \Delta_n} \varphi_{(\beta,j)} \gamma_{(\beta,j)}
\end{equation}
with $\varphi_{(\beta,j)}\in \ef $. \\
Again, we have the isomorphism :
\begin{align}\label{isom2}
\begin{split}
\Phi: \ef^{\Delta_n}&\longleftrightarrow \Hom_{\ef}(\Hom_{\ef}(G_n,\ef),\ef)\\
  y=(y_{(\beta,j)})_{(\beta,j)_\in\Delta_n}&\longmapsto\left\{
                                                        \begin{array}{ll}
                                                         \Phi(y):\Hom_{\ef}(G_n,&\ef)\rightarrow \ef \\
                                                           &\gamma_{(\beta,j)} \mapsto \Phi(y)(\gamma_{(\beta,j)})=y_{(\beta,j)}.
                                                        \end{array}
                                                      \right.
\end{split}
\end{align}
Consider the restrictions
\begin{align}\label{lambda-n}
\begin{split}
w_\lambda|_{G_n}&=w_\lambda^{(n)}\quad\text{with}\;\langle \delta_\alpha e_i,w_{\lambda}^{(n)}\rangle=w_{\lambda (\alpha,i)},\\
q|_{G_n}&=q_n\quad\text{with}\;\langle\delta_\alpha e_i,q_n\rangle=q_{(\alpha,i)}
\end{split}
\end{align}
for $n\in\enn, \gamma\in\Lambda$ and ${(\alpha,i)}\in\Delta_n$. We will need the following lemma :
\begin{lem}\label{lemme}
$q_n\in \langle (w_\lambda^{(n)})_{\lambda\in \Lambda}\rangle$.
\end{lem}
\prf We have $q_n\in \Hom_{\ef}(G_n,\ef)$; suppose that $q_n\notin\langle (w_\lambda^{(n)})_{\lambda\in \Lambda}\rangle$. In this case, we know that there exists   $\Theta\in  \Hom_{\ef}(\Hom_{\ef}(G_n,\ef),\ef)$ with $\Theta(w_\lambda^{(n)})=0$ for $\lambda\in \Lambda$ and $\Theta(q_n)=1$. By the isomorphism  \eqref {isom2}, there exists $y=(y_{(\alpha,i)})_{(\alpha,i)\in\Delta_n}\in \ef^{\Delta_n}$ such that $\Phi(y)=\Theta$, i.e
$\Theta(\gamma_{(\alpha,i)})=y_{(\alpha,i)}$ for  $(\alpha,i)\in\Delta_n\ $. Take
\begin{equation*}
g=\sum_{\alpha\in \Delta_n }y_{(\alpha,i)} \delta_\alpha e_i\ \in G_n.
\end{equation*}
 Then
\begin{align}
\begin{split}
    \langle g,w_\lambda \rangle&=\sum_{(\alpha,i)\in\Delta_n}y_{(\alpha,i)} w^{(n)}_{\lambda(\alpha,i)}\\
    &=\sum_{(\alpha,i)\in\Delta_n}\Theta(\gamma_{(\alpha,i)})w^{(n)}_{\lambda(\alpha,i)}\\
    &=\sum_{(\alpha,i)\in \Delta_n}\Theta(w^{(n)}_{\lambda(\alpha,i)}\gamma_{(\alpha,i)})\\
    &=\Theta(\sum_{(\alpha,i)\in  \Delta_n}w^{(n)}_{\lambda(\alpha,i)}\gamma_(\alpha,i))\\
    &=\Theta(w_\lambda^{(n)})=0.
\end{split}
\end{align}

We then get
\begin{equation*}
(\forall \lambda\in\Lambda) \quad \bigg[\langle g,w_\lambda\rangle=0\bigg],
\end{equation*}
which proves that $g\in V^{\perp}$. But we also have
\begin{align*}
\langle g,q\rangle&=\sum_{(\alpha,i)\in\Delta_n}y_{(\alpha,i)} q_{(\alpha,i)}\\
&= \sum_{{(\alpha,i)}\in\Delta_n}\Theta (\gamma_{(\alpha,i)}) q_{(\alpha,i)}=\Theta(\sum_{{(\alpha,i)}\in\Delta_n}\gamma_{(\alpha,i)} q_{(\alpha,i)}),
\end{align*}
and by \eqref{phi},
\begin{equation*}
  q_n=\sum_{(\alpha,i)\in\Delta_n}\gamma_{(\alpha,i)} q_{n_{(\alpha,i)}}=\sum_{(\alpha,i)\in\Delta_n}\gamma_{(\alpha,i)} q_{(\alpha,i)},
\end{equation*}
so that
\begin{equation*}
\langle g,q\rangle =\Theta(q_n)=1.
\end{equation*}
 Since $q\in V^{\perp\perp}$, this implies that $g\notin V^{\perp\perp\perp}=V^\perp$, which is a contradiction. We conclude that necessarily   $q_n\in \langle (w_\lambda^{(n)})_{\lambda\in \Lambda}\rangle$.\qed\\

\noindent\textit{Proof of theorem \ref{closed}} (continued). For  $n \in \enn^*$, there is then a family $(\mu _{\lambda}^{(n)})_{\lambda\in \Lambda}$ with finite support such that
 \begin{equation*}
    q|_{G_n}=q_n=\sum_{\lambda\in \Lambda} \mu _{\lambda}^{(n)} w_\lambda^{(n)}.
\end{equation*}
Consider  the element $q'_n=\sum_{\lambda\in \Lambda} \mu _{\lambda}^{(n)} w_\lambda\in V$, which is an extension of $q_n$ on $\aaa^l$. For $(\alpha,i)\in \Delta_n$, we have, by \eqref{lambda-n}
 \begin{equation*}
  q_{(\alpha,i)}={q'_n}_{(\alpha,i)},
\end{equation*}
so that $q-q'_n\in V_n$. Therefore,  the sequence  $(q'_n)_{n\in \enn^*}$ converges to $q$ in $\aaa^l$. Finally, we have  $q\in \overline{V}$. Hence $V^{{\perp}{\perp}}\subset V$ and the equality holds. Taking $G=V^{\perp}$, we get $G^{\perp}=V^{{\perp}{\perp}}=V$ and the theorem is proved.\qed

\section{Use of polynomials  and power series }\label{section5}

Let $r\geqslant1$  be an integer. For $\rho=1\ldots,r,$ let $X_\rho$ (resp. $Y_\rho$) be letters (or variables). For $\alpha\in\enn^r$ we define $X^\alpha$ (resp. $Y^\alpha$) by
\begin{equation}\label{X-Y-alpha}
  X^\alpha=X_1^{\alpha_1}\cdots X_r^{\alpha_r}\quad\text{(resp. } Y^\alpha=Y_1^{\alpha_1}\cdots Y_r^{\alpha_r}\text{)}.
\end{equation}

Let $\df=\ef[X_1,\ldots ,X_r]$  be the $\ef$-vector space of the polynomials with the $r$ variables $X_1,\ldots, X_r$ and $\af=\ef[[Y_1,\ldots, Y_r]]$ the  $\ef$-vector space of  the formal power series  with the $r$ variables $Y_1,\ldots, Y_r$. The family $(X^\alpha)_{\alpha\in\enn^r}$ is  $\ef$-base of $\df$, thus an element of $\df$ can be written uniquely as
\begin{equation*}
d(X)=\sum_{\alpha\in\enn^r} d_\alpha X^\alpha\quad  \text{ with}\quad  d_\alpha\in\ef\quad  \text{for all } \alpha\in\enn^r,
\end{equation*}
where $d_\alpha=0$ except for a finite number of $\alpha$. An element $W(Y)$ of $\af$ can be expressed uniquely as
  \begin{equation*}
W(Y)=\sum_{\alpha\in\enn^r} W_\alpha Y^\alpha\quad\text{ with}\quad d_\alpha\in\ef\quad\text{for all}\quad \alpha\in\enn^r.
 \end{equation*}
Therefore, {using } the  sets $\aaa$ and $\ddd$ introduced in Section \ref{section2} , we obtain  the following vector spaces isomorphisms
\begin{align}\label{1-d-df-a-af-isom}
\begin{split}
\mathcal{D}: \df=\ef[X_1,\ldots,X_r]&\cong \ddd\\
d(X)&=\sum_{\alpha\in\enn^r} d_\alpha X^\alpha\longleftrightarrow(d_{\alpha})_{\alpha\in\enn^r}, \\
\mathcal{A} : \af=\ef[[Y_1,\ldots, Y_r]]&\cong \aaa\\
W(Y)=\sum_{\alpha\in\enn^r} W_\alpha Y^\alpha&\longleftrightarrow W=(W_\alpha)_{\alpha\in\enn^r}.
\end{split}
\end{align}
By these isomorphisms, we may identify $X^\alpha$ (resp. $Y^\alpha$) with the element $\delta_\alpha$ of $\ddd$ (resp. of $\aaa$). If $d(X)\in\df^l$ with $d(X)=\sum_{\alpha\in\enn^r}d_\alpha X^\alpha$ for $\lambda=1,\ldots,l$, we then may write the equalities
\begin{equation}\label{pol-ident}
d(X)=\sum_{\alpha\in\enn^r} d_\alpha X^\alpha=(d_\alpha)_{\alpha\in\enn^r}, \end{equation}
  and for  $W(Y)=\sum_{\alpha\in\enn^r}W_\alpha Y^\alpha\in\af$, we  may write the equalities
\begin{equation}\label{power_ident}
W(Y)=W=(W_{\alpha})_{\alpha\in\enn^r}.
\end{equation}
Using \eqref{d-df-a-af-isom} and proposition \ref{scal-prod}, we obtain the scalar product
\begin{align}
\begin{split}\label{a-d-scalar-prod}
\langle -,-\rangle :\df\times\af &\longrightarrow\ef\\
(d(X),W(Y))&\longmapsto\langle d(X),W(Y)\rangle=d(X)\cdot W(Y)= \sum_{\alpha\in\enn^r} d_\alpha W_\alpha .
\end{split}
\end{align}
and  the isomorphism
\begin{equation}\label{a-d-iso}
\af\cong\Homf(\df,\ef),\ W(Y)\longmapsto\langle-,W(Y)\rangle_{\ef}.
\end{equation}
Considering $W\in \af$ as element of $\Homf(\df,\ef)$, we have
\begin{equation}\label{W-X-alpha-1}
\langle d(X),W\rangle=W(d(X)).
\end{equation}
In other terms, with the identification $\af=\Homf(\df,\ef)$, we may view $W$ as acting on $X^\alpha$ by
\begin{equation}\label{W-X-alpha-2}
W(X^\alpha)=W(\alpha)=W_\alpha.
\end{equation}

For   an integer $l\geqslant 1$, the set $\df^l$ is those  of the  polynomial vector with $l$ columns, with entries in $\df$ and $\af^l$ is those  of  the column-vector formal powers series with $l$ rows, with entries in $\af$. An element of $\df^l$ is of the form
\begin{equation}\label{elt-df}
  d(X)=(d_1(X),\ldots,d_l(X))\quad\text{with}\quad d_\lambda(X)=\sum_{\alpha\in \enn^r}^ld_{\lambda\alpha}X^\alpha\in\df,
\end{equation}
for $\lambda=1,\ldots,l$ and an element of $\af^l$ is of the form
\begin{equation}\label{elt-df}
 W(Y)=\begin{pmatrix}W_1(Y)\\
\vdots \\
W_l(Y)\\
\end{pmatrix}\quad\text{with}\quad W_\lambda(Y)=\sum_{\alpha\in\enn^r}W_{\lambda\alpha}Y^\alpha\in\af
\end{equation}
for $\lambda=1,\ldots,l$. \\

We have the isomorphism of vector spaces
\begin{align}\label{l-d-df-a-af-isom}
\begin{split}
\mathcal{D}: \df^l&\cong \ddd^l\\
d(X)=(d_1(X),\ldots,d_l(X))&\longleftrightarrow((d_{1\alpha})_{\alpha\in\enn^r},\ldots,(d_{l\alpha})_{\alpha\in\enn^r}) \\
\mathcal{A} : \af^l&\cong \aaa^l\\
W(Y)=\begin{pmatrix}W_1(Y)\\
\vdots \\
W_l(Y)
\end{pmatrix}&\longleftrightarrow W=\begin{pmatrix}(W_{1\alpha})_{\alpha\in\enn^r}\\
\vdots \\
(W_{l\alpha})_{\alpha\in\enn^r}
\end{pmatrix}..
\end{split}
\end{align}
Note that $\df^l$ is a $\df$-module, with the multiplication

\begin{equation}\label{D-module}
d(X)\cdot (d_1(X),\ldots,d_l(X))=(d(X)d_1(X),\ldots,d(X)_ld(X)),
\end{equation}
for $d(X)\in\df^l$ and $(d_1(X),\ldots,d_l(X))\in\df^l$. \\

By \eqref{l-scalar-prod} and \eqref{d-df-a-af-isom}, we obtain the scalar product
\begin{align}\label{l-a-d-scalar-prod}
  \begin{split}
  \langle-,-\rangle:\df^l\times\af^l&\longrightarrow\ef,\\
  (d(X),W(Y))&\longmapsto \langle d(X),W(Y)\rangle = \sum_{\lambda=1}^{l}(\sum_{\alpha\in\enn^{r}}d_{\lambda\alpha}\cdot W_{\lambda\alpha})
 \end{split}
\end{align}
and then the vector spaces isomorphism
\begin{align}\label{l-af-df-isom}
 \begin{split}
\af^l&\cong\Homf(\df^l,\ef)\\
 W(Y)=\begin{pmatrix}W_1(Y)\\
 \vdots \\
W_l(Y)
\end{pmatrix} &\longmapsto\langle-,W(Y)\rangle.
 \end{split}
\end{align}
For $P\subset\af^l$ and $Q\subset\df^l$, we define the orthogonals $P^{\perp}$ and $Q^{\perp}$ by
\begin{align}
  P^{\perp} &=\{d(X)\in \df^l\ |\langle d(X),W(Y)\rangle=0\; \forall W(Y)\in P\} \subset \df^l\label{P-l-orth}, \\
  Q^{\perp} &=\{W(Y)\in \af^l\ |\langle d(X),W(Y)\rangle\ \forall\; d(X)\in Q\} \subset\af^l\label{Q-l-orth}.
\end{align}
Let $P\subset \df^l$ a non-empty subset, $\mathcal{D}$ and $\mathcal{A}$ the vector spaces isomorphisms in \ref{l-d-df-a-af-isom}. Set $P_{0}=\mathcal{D}(P)\subset\ddd^l$, with
\begin{equation*}
  P_0=\Biggl\{((d_{1\alpha})_{\alpha\in\enn^r},\ldots,(d_{l\alpha})_{\alpha\in\enn^r})\;\Bigg\vert\;(\sum_{\alpha\in\enn^r}d_{1\alpha}X^\alpha,\ldots,\sum_{\alpha\in\enn^r}d_{l\alpha}X^\alpha\rangle\in P\Biggr\}.
\end{equation*}
Then the mapping
\begin{align*}
 P&\longleftrightarrow P_0\\
 (d_1(X),\ldots,d_l(X))&\longleftrightarrow((d_{1\alpha})_{\alpha\in\enn^r},\ldots,(d_{l\alpha})_{\alpha\in\enn^r})
\end{align*}
is an isomorphism of vector spaces.
\begin{prop}\label{A-P}With the above notations, one has $P^\perp_0=\mathcal{A}(P^\perp)$, which leads to the following vector spaces isomorphism
\begin{align*}
 P^\perp&\longleftrightarrow P_0^\perp=\mathcal{A}(P^\perp)\\
 W(Y)=\begin{pmatrix}W_1(Y)\\
\vdots \\
W_l(Y)
\end{pmatrix}&\longleftrightarrow W=\begin{pmatrix}(W_{1\alpha})_{\alpha\in\enn^r}\\
\vdots \\
(W_{l\alpha})_{\alpha\in\enn^r}
\end{pmatrix}.
\end{align*}
\end{prop}
\prf It is left to the reader.\qed\\

Now, we are going to introduce a topology in $\af^l$.

\begin{defi}[Convergence in  $\af^l$] A sequence
$(W_n(Y))_{n\in\enn}$ of elements of $\af^l$ converges to the element $W(Y)\in\af^l$ if the sequence of elements $(W_n)_{n\in\enn}$ of  $\aaa^l$, where
$$W_n= (W_{n\alpha})_{\alpha\in\enn^r}$$
converges to the element $W=(W_{\alpha})_{\alpha\in\enn^r}$ of $\aaa^l$.
\end{defi}
\begin{defi}[Topology on $\af^l$]\label{l-A-top}The  topology of $\af^l$, which is the topology of the \textit{pointwise convergence} is defined as follows: a subset $F$ of $\af^l$  is closed if $F=\emptyset$ or $F=\af^l$ or it verifies the following property : if a sequence  $(W_n(Y))_{n\in\enn}$ of elements of $F$ converges to $W(Y)\in\af^l$, then $W(Y)\in F$.
\end{defi}
One can verify that the sets $F$'s in Definition \ref{l-A-top} verify the axioms of  closed spaces (closure under finite union and arbitrary intersection), so that one can indeed define a topology whose closed sets are these sets.

\begin{thm}\label{orth-closed} A subset $V$ of $\af^l$ is closed if and only if it is of the form  $P^\perp$ where $P$ is a subset of $\df^l$.\end{thm}

\prf Suppose that $V\neq\emptyset$ is closed in $\af^l$. Then $V_0=\mathcal{A}(V)$ is closed in $\aaa^l$. We know, by \ref{l-closed} that $V_0$ is of the form $G_0^\perp$, where $G_0\subset\ddd^l$. By \ref{A-P}, there exists $P\subset\df^l$ such that $G_0=\mathcal{D}(G)$, so that $V_0=G_0^\perp=\mathcal{A}(G^\perp)$. We conclude that $V=G^\perp$.\\
Conversely, let $P$ be a subset of $\df^l$. If $P=\emptyset$, then $P^\perp=\emptyset$ is closed in $\af^l$. Suppose that $P\neq\emptyset$.  Let $(W_n(Y))_{n\in\enn}$ a sequence of elements of $P^\perp$ which converges to $W(Y)\in\af^l$. Using  the notations in \ref{A-P}, this means that, the sequence
$(W_n)_{n\in\enn}$ of $P_0^\perp$ converges to $W$ in $\aaa^l$. Since $P_0^\perp$ is closed in $\aaa^l$, it follows that $W\in P_0^\perp$.  Hence, $W(Y)\in P^\perp$. We have shown that $P^\perp$ contains all the limits of all its convergent sequences, i.e. $P^\perp$  is closed. \qed\\

\section{Shift operator, polynomial operator in the shift and polynomial-power series multiplication}\label{section6}
Let $\vect(\ef)$ the category whose object consists of all $\ef$-vector spaces and for two objects $E,F\in\vect(\ef)$, the set of morphism from $E$ to $F$ is $\Homf(E,F)$, which consists of all linear mappings from $E$ to $F$. We then have the \textit{covariant functor}  $\Homf(-,\ef)$ defined by
\begin{align}\label{adj}
\begin{split}
\Homf(-,\ef) :\vect(\ef)&\longrightarrow\vect(\ef)\\
E&\longmapsto\Homf(E,\ef)\\
(f :E\longrightarrow F)&\longmapsto\left\{
                                     \begin{array}{ll}
                                      \Homf(f,\ef) :\Homf(F,&\ef)\longrightarrow \Homf(E,\ef) \\
                                        &u\longmapsto u\circ f,
                                     \end{array}
                                   \right.
\end{split}
\end{align}
(\cite{la}). This leads to the following definition :
\begin{defi}[\cite{la}]Let $E,F\in\vect(\ef)$ and $f\in\Homf(E,F)$. The \textit{functorial adjoint} of $f$ if the linear mapping $\Homf(f,\ef)$.
\end{defi}
Now we are going to look the adjoints of particular linear mappings: take $E=F=\df$ and take $d(X)\in\df$. We get the \textit{multiplication by} $d(X)$, which is the linear mapping
\begin{equation*}
d(X): \df\longrightarrow \df,\; c(X)\longmapsto c(X)\cdot d(X)
\end{equation*}
which we also denoted by $d(X)$. For the case  $d(X)=X^\beta$ and $\beta\in\enn^r$. We get the `` multiplication by $X^\beta$ '':
\begin{equation*}
X^\beta:\df\longrightarrow\df,\; c(X)\longmapsto c(X)\cdot X^\beta.
\end{equation*}
The adjoint of the multiplication by $X^\beta$ is given by the following lemma:
\begin{lem}
The functorial adjoint  of the multiplication by $X^\beta$
\begin{align*}
\begin{split}
X^\beta:\df&\longrightarrow\df\\
c(X)&\longmapsto c(X)\cdot X^\beta,
\end{split}
\end{align*}
is the  $\ef$-endomorphism
\begin{align}\label{decal1}
\begin{split}
 \af&\longrightarrow \af\\
 W(Y)=\sum_{\alpha\in\enn^r} W_\alpha Y^\alpha&\longmapsto\sum_{\alpha\in\enn^r} W_{\alpha+\beta}Y^\alpha.
 \end{split}
\end{align}
\end{lem}
\prf
Using \eqref{a-d-iso}, we obtain $\Homf(\df,\ef)=\af$. Since for  $W\in\af$, the map $\Homf(X^\beta,\ef)(W) =W\circ X^\beta$ is an element of $\af$
and
\begin{equation}\label{muladj}
(W\circ X^\beta)(\alpha)=W\circ X^\beta(X^\alpha)=W(X^\beta\cdot X^\alpha)=W(X^{\alpha+\beta})=W_{\alpha+\beta}
\end{equation}
for $\alpha\in\enn^r$ (see \eqref{W-X-alpha-1} and \eqref{W-X-alpha-2}), this completes the proof of  \eqref{decal1}. \qed

The adjoint of the multiplication by $X^\beta$ is the \textit{shift operation}, \cite{ob,wil1,wil2,wil3,wil5}, and also denoted by $X^\beta$. We use symbol ``$\circ$ '' to mean that actually, $X^\beta$ operates on a power series. Thus,
\begin{equation*}
X^\beta\circ W(Y)=\Homf(X^\beta,\ef)(W(Y))= \sum_{\alpha\in\enn^r} W_{\alpha+\beta} Y^\alpha \in\af.
\end{equation*}
\begin{ex}Fix $\alpha\in\enn^r$ and take $W(Y)=Y^\alpha$. For $\beta\in\enn^r$, we have
\begin{equation}\label{exa-shift}
X^\beta\circ Y^\alpha=\left\{
                        \begin{array}{ll}
                          Y^{\alpha-\beta} & \hbox{if}\quad\beta\leqslant_+\alpha, \\
                          0 & \hbox{otherwise},
                        \end{array}
                      \right.
\end{equation}
where $\beta\leqslant_+\alpha$ means that $\beta_i\leqslant\alpha_i$ for all $i=1,\ldots,r$.
\end{ex}
We have the following fundamental property:
\begin{equation}\label{shift-prop}
    (X^\alpha \cdot X^\beta)\circ W(Y)=X^\alpha \circ(X^\beta\circ W(Y))
\end{equation}
for all $\alpha,\beta\in \enn^r$ and $W(Y)\in \af$.\\

Now consider the general case of the polynomial multiplication by $d(X)\in\df$:
\begin{equation*}
    d(X):\df\longrightarrow\df,\ \ c(X)\longmapsto c(X)\cdot d(X).
\end{equation*}
If $d(X)=\sum_\beta d_\beta X^\beta$, we may view $d(X)$ as a linear combination of  $X^\beta$. Taking the adjoint, we have
\begin{equation*}
\Homf(d(X),\ef)=\sum_\beta d_\beta \Homf(X^\beta ,\ef),
\end{equation*}
and using \eqref{decal1}, we have
\begin{align}
\begin{split}
    \Homf(d(X),\ef)(W(Y))&=\sum_\beta d_\beta\Homf(X^\beta,\ef)(W(Y))\\
&=\sum_\alpha(\sum_\beta(d_\beta W_{\alpha+\beta})Y^\alpha,
\end{split}
\end{align}

We may view this last equation as a generalisation of the shift operation and consider it as an operation of $d(X)$ on $W(Y)$ : this is \textit{the polynomial operator in the shift}, \cite{ob,wil1,wil2,wil3,wil5}. We have proved the following theorem :
\begin{thm}\label{fundamental} The functorial adjoint of the polynomial multiplication by $d(X)$
\begin{align}
\begin{split}
    d(X):\df&\longrightarrow\df\\
    c(X)&\longmapsto c(X)\cdot d(X),
\end{split}
\end{align}
is the\textit{ polynomial operator in the shift}, also denoted by $d(X)$ and defined as
\begin{align}\label{po-oper-sh}
\begin{split}
d(X) : \af &\longrightarrow \af \\
 W(Y)&\longmapsto d(X)\circ W(Y)=\sum_\alpha(\sum_\beta d_\beta W_{\alpha+\beta})Y^\alpha.
\end{split}
\end{align}
\end{thm}
This leads to an operation of $\df$ on $\af$, given as follows :
\begin{defi}\label{D-A-op}
The polynomial-power series \textit{multiplication }$\circ$ is the operation
\begin{align}
\begin{split}
 \df \times \af &\longrightarrow \af \\
 (d(X),W(X)) &\longmapsto d(X)\circ W(Y)
\end{split}
\end{align}
where
\begin{equation*}
 d(X)\circ W(Y)=\Homf(d(X),\ef)(W(Y))=\sum_\alpha(\sum_\beta d_\beta W_{\alpha+\beta})Y^\alpha.
\end{equation*}
\end{defi}
Using the fundamental property \eqref{shift-prop}, $\circ$ is indeed an external operation of $\df$ on $\af$. Moreover, its has a more interesting properties, whose proofs are left to the reader :
\begin{prop}The multiplication $\circ$ provides $\af$ and $\af^l$, for all integer $l\geqslant 1$, with a $\df$-module structure.
\end{prop}

In \cite{b13,b1}, a generalization of the definition \ref{D-A-op} by taking matrices of polynomials and powers series  is given. For this, we introduce the following notations: for integers $k,l\geqslant 1$, the
   $\df^{k,l}$ is those of the polynomial matrices with $k$ rows and $l$ columns, with coefficient in $\df$.  An element of $\df^{k,l}$ is of the form
\begin{equation}\label{mat-pol}
  R(X)=(R_{\kappa\lambda}(X))_{\substack{1\leqslant\kappa\leqslant k\\
  1\leqslant\lambda \leqslant l}}
\end{equation}
with
\begin{equation}\label{mat-R}
R_{\kappa\lambda}(X)=\sum_{\alpha\in\enn^r}R_{\kappa\lambda\alpha}X^\alpha\in\df
\end{equation}
for $\kappa=1,\ldots,k$ and $\lambda=1,\ldots,l$.
   For simplicity, the set of row vectors with $l$ columns $\df^{1,l}$ is simply denoted by $\df^l$.\\
\begin{thm}[\cite{b13,b1,ob}]\label{l-fundamental}Fix   $R(X)=(R_{\kappa\lambda}(X))_{\substack{1\leqslant\kappa\leqslant k\\
  1\leqslant\lambda \leqslant l}}\in\df^{k,l}$ where $R_{\kappa\lambda}(X)\in\df$ is as in \ref{mat-R}.  The functorial adjoint of the polynomial-polynomial matrix multiplication
\begin{align}\label{polmult}
\begin{split}
R(X)^T : \df^k &\longrightarrow \df^l \\
c(X)&\longmapsto c(X)\cdot R(X),
\end{split}
\end{align}
denoted by $R(X)$ is the $\df$-linear mapping
\begin{align}\label{mpol-1}
\begin{split}
R(X) : \af^l &\longrightarrow \af^k \\
 W(Y)&\longmapsto R(X)\circ W(Y)=\left(
           \begin{array}{c}
             \sum_{\lambda=1}^l\sum_{\beta}(\sum_{\alpha} R_{1\lambda\alpha}W_{\lambda(\beta+\alpha)}) Y^\beta \\
             \vdots \\
             \sum_{\lambda=1}^l\sum_{\beta}(\sum_{\alpha} R_{\kappa\lambda\alpha}W_{\lambda(\beta+\alpha)}) Y^\beta \\
             \vdots\\
             \sum_{\lambda=1}^l\sum_{\beta}(\sum_{\alpha} R_{k\lambda\alpha}W_{\lambda(\beta+\alpha)}) Y^\beta
           \end{array}
         \right).
\end{split}
\end{align}
The mapping $R(X)$ is also the algebraic adjoint of the mapping $R(X)^T$. This means that the $\ef$-linear mapping \eqref{adj} is the unique linear mapping from  $\aaa^l$ to $\aaa^k$ which verifies
\begin{equation*}
  \langle d(X)R(X),W(Y)\rangle = \langle d(X), R(X)\circ W(Y)\rangle
\end{equation*}
\end{thm}
By \eqref{mpol-1}, the $\kappa$-th row of the vector of power series $R(X)\circ W(Y)$ is
\begin{equation}\label{kappa-polmult}
  [R(X)\circ W(Y)]_\kappa=\sum_{\beta}(\sum_{\alpha}\sum_{\lambda=1}^l R_{\kappa\lambda\alpha}W_{\lambda(\beta+\alpha)}) Y^\beta
\end{equation}
and the coefficient of $Y^\beta$ is
\begin{equation}\label{l-kappa-coeff-Y}
   [R(X)\circ W(Y)]_{\kappa_\beta}=\sum_{\alpha}\sum_{\lambda=1}^l R_{\kappa\lambda\alpha}W_{\lambda(\beta+\alpha)}.
\end{equation}
The following corollary is immediate :
\begin{cor}The mapping
\begin{align}\label{d-a-bilin}
\begin{split}
    \circ:\df^{k,l}\times \af^l&\longrightarrow\af^k\\
 (R(X),W(Y))&\longmapsto R(X)\circ W(Y).
\end{split}
\end{align}
is bilinear.
\end{cor}

\noindent Setting $k=1$ in the above corollary, we get the $\df$-bilinear mapping
\begin{align}\label{A-scal-prod}
\begin{split}
\df^l\times\af^l&\longrightarrow\af\\
(d(X),W(Y)) &\longmapsto \langle d(X),W(Y) \rangle_{\af}=d(X)\circ W(Y)
\end{split}
\end{align}
 which is also a scalar product (\cite{ob}). Thus, we can also
define orthogonals with respect to this bilinear mapping   which we  denote by ${\perp_{\af}}$ : for a subset $P$ (resp. $Q$) in $\df^l$ (resp.$\af^l$), we have
\begin{align}\label{orthrond}
\begin{split}
&P^{\perp_{\af}} =\{W(Y)\in\af^l|\;p(X)\circ W(Y)=0\; \forall p(X)\in P\}\subset\af^l,\\
 &Q^{\perp_{\af}} = \{p(X)\in\df^l|\;p(X)\circ W(Y)=0\; \forall W(Y)\in Q\}\subset\df^l.
\end{split}
\end{align}
The proof of the following remark is left to the reader:
\begin{rem}\label{orth-sm}(1) Let $\langle P\rangle $ (resp. $\langle Q\rangle$) the $\df$-submodule of $\df^l$ (resp. of $\af^l $ generated by $P$ (resp. by $Q$). Then
\begin{equation}
  P^{\perp_{\af}}=\langle P\rangle^{\perp_{\af}}\quad (\text{resp. }Q^{\perp_{\af}}=\langle Q\rangle^{\perp_{\af}} ).
\end{equation}
Note that $P^{\perp_{\af}}$ (resp. $Q^{\perp_{\af}}$) is an $\df$-submodule of $\af^l$ (resp. of $\df^l$).\\
(2) Given $R(X)\in\df^l$  with $R(X)=(R_\lambda(X))_{\lambda=1,\ldots,l}$ (a polynomial vector; where the index $\kappa$ is omitted in \eqref{mat-pol}), according to \eqref{l-kappa-coeff-Y} the coefficient of $Y^\beta$ is
\begin{equation}\label{kappa-coeff-Y}
   [R(X)\circ W(Y)]_\beta=\sum_{\alpha}\sum_{\lambda=1}^l R_{\lambda\alpha}W_{\lambda(\beta+\alpha)}.
\end{equation}
(3) For $d(X)\in\df^l$ and $W(Y)\in\af^l$, comparing \eqref{l-a-d-scalar-prod} with \eqref{kappa-coeff-Y} and \eqref{A-scal-prod}  yields
\begin{equation}\label{A-scal-prod-coeff}
  \langle d(X),W(Y)\rangle=[d(X)\circ W(Y)]_0=\langle d(X),W(Y)\rangle_{\af}.
\end{equation}
(4) If $P\subset\df^l$, then the above equation \eqref{A-scal-prod-coeff} yields
\begin{equation}\label{orths-incl}
 P^{\perp_{\af}}\subset P^\perp.
\end{equation}
\end{rem}
The following proposition gives a condition for the converse of \eqref{orths-incl}:
\begin{prop}\label{recipr-orths}If $P\subset\df^l$ and $P^\perp$ is a $\df$-submodule of $\af^l$, then $P^\perp=P^{\perp_{\af}}$.
\end{prop}
\prf We know from \eqref{orths-incl} that $P^{\perp_{\af}}\subset P^\perp$. Now, fix $W(Y)\in P^\perp$.
Since $P^\perp$ is a $\df$-submodule of $\af^l$, it also verifies $X^\gamma\circ W(Y)\in P^\perp$ for $\gamma\in \enn^r$, so that
for $d(X)\in P$, we have
\begin{equation}\label{d-X-gamma}
 d(X)\cdot (X^\gamma \circ W(Y))=0.
\end{equation}
Write $d(X)=(d_1(X),\ldots,d_l(X))$ with $d_\lambda(X)=\sum_{\alpha}d_{\lambda\alpha} X^\alpha$ for $\lambda=1,\ldots,l$ and
\begin{equation*}
  W(Y)=\begin{pmatrix}W_1(Y) \\
\vdots \\
W_\lambda(Y)\\
\vdots \\
W_l(Y) \\
\end{pmatrix}\quad\text{with}\quad W_\lambda(Y)=\sum_{\alpha\in\enn^r}W_{\lambda\alpha}Y^\alpha.
\end{equation*}
Rewriting \eqref{d-X-gamma} yields
\begin{align*}
 0=d(X)\cdot (X^\gamma \circ W(Y))&=\langle d(X), X^\gamma\circ W(Y)\rangle\\
 &=\langle (d_1(X),\ldots,d_l(X)),\begin{pmatrix} \\
\sum_{\alpha\in\enn^r}W_{1(\alpha+\gamma)}Y^\alpha \\
\vdots \\
\sum_{\alpha\in\enn^r}W_{l(\alpha+\gamma)}Y^\alpha \\
\end{pmatrix}\rangle\\
&=\sum_{\alpha\in\enn^r}\sum_{\lambda=1}^l d_{\lambda\alpha}W_{\lambda(\gamma+\alpha)}.
\end{align*}
The comparison with   \eqref{l-kappa-coeff-Y}, where $R(X)$ is set to $d(X)$ ($k=1$ and $\kappa$ is omitted)
yields
\begin{equation*}
0=d(X)\cdot (X^\gamma \circ W(Y))=[d(X)\circ W(Y)]_\gamma.
\end{equation*}
Since it is true for $\gamma\in\enn^r$, we have $d(X)\circ W(Y)=0$  for $d(X)\in P$, so that $W(Y)\in P^{\perp_{\af}}$. It follows that $P^\perp\subset P^{\perp_{\af}}$ and the equality holds.\qed
\begin{prop}Given $R(X)\in\df^{k,l}$, the kernel of the $\df$-linear mapping
\begin{align}
\begin{split}
    R(X): \af^l&\longrightarrow\af^k\\
W(Y)&\longmapsto R(X)\circ W(Y)
\end{split}
\end{align}
is the $\df$-module
 \begin{equation}\label{ker-R}
 \Ker R(X)=\{W(Y)\in\af^l|R(X)\circ W(Y)=0\}\subset\af^l.
 \end{equation}
 \end{prop}
These modules were given a special name.
\begin{defi}
A $\df$-submodule of $\af^l$ of the form $\Ker R(X)$ where $R(X)\in\df^{k,l}$ is called \textit{autoregressive module}.
\end{defi}

The next proposition gives characterizes the autoregressive modules.
\begin{prop}\label{auto}For  a nonempty subset $P$ of $\df^l$, the module $P^{\perp_{\af}}$ is autoregressive. Conversely, an autoregressive module is of the form $P^{\perp_{\af}}$,  were $P\subset\df^l$ for some nonnegative integer $l$.
\end{prop}
\prf According to the Hilbert basis theorem  (\cite{b11}), the submodule $\langle P\rangle$ of $\df^l$ is finitely generated: there exists $k\in\enn^*$ and a polynomials $R_1(X),\ldots,R_k(X)\in \df^l$ (with $ R_\kappa(X)=(R_{\kappa\lambda}(X))_{1\leqslant\lambda\leqslant l}$ for  $ \kappa=1,\ldots,k$) such that $\langle P\rangle$ is generated by these polynomials:
\begin{equation*}
\langle P\rangle=\langle R_1(X),\ldots,R_k(X)\rangle.
\end{equation*}
It follows that
\begin{equation*}
 \langle P \rangle^{\perp_{\af}}=\langle R_1(X),\ldots,R_k(X)\rangle^{\perp_{\af}}.
\end{equation*}
But, according to \eqref{orth-sm}, we get
\begin{equation*}
    P^{\perp_{\af}}=\langle P\rangle^{\perp_{\af}}=\langle R_1(X),\ldots,R_k(X)\rangle^{\perp_{\af}}.
\end{equation*}
Finally, we have
\begin{eqnarray*}
  P^{\perp_{\af}} &=& \{W(Y)\in\af^l|R_\kappa(X)\circ W(Y)=0\ \forall \kappa=1,\ldots,k \} \\
    &=&  \{W(Y)\in\af^l|R(X)\circ W(Y)=0\}\\
    &=&\Ker R(X)\quad\text{(see \eqref{ker-R})}.
\end{eqnarray*}
Conversely if $M=\Ker R(X)$ with $R(X)\in\df^{k,l}$ is an autoregressive module, then the above proof shows that $M=P^{\perp_{\af}}$, where $P$ is the set of the rows of  $R(X): P=\{ R_1(X),\ldots,R_k(X)\}$
. \qed\\

\section{Discrete dynamical systems with $r$ shifts}\label{section7}
We present here a synthesis of the various definitions of discrete algebraic systems given in \cite{ob,wil1,wil2,wil3,wil5}.

\begin{thm}\label{autoregressive-orth-closed}
Let $\boldsymbol{\mathcal{B}}$ be a subset of $\af^l$. Then the following properties are equivalent:\\
$(1)$ There exists $R(X)\in\df^{k,l}$ such that $\boldsymbol{\mathcal{B}}=\Ker R(X)$.\\
$(2)$ There exists a non-empty finite subset $P$ of $\df^l$ such that $\boldsymbol{\mathcal{B}}=P^{\perp_{\af}}$.\\
$(3)$ The set $\boldsymbol{\mathcal{B}}$ is a $\df$-submodule of $\af^l$ which is closed for its topology.
\end{thm}
\prf $(1)\Longleftrightarrow (2)$ : this is Proposition \ref{auto}.\\
$(2)\Longrightarrow (3)$: If $\boldsymbol{\mathcal{B}}=P^{\perp_{\af}}$ with $P\subset\df^l$, we know by  (1), Remark \ref{orth-sm} that $P^{\perp_{\af}}$ is a $\df$-submodule of $\af^l$. We will show that it is close in $\af^l$. Let $(W_n(Y))_{n\in\enn} $ be a sequence  of elements of $\boldsymbol{\mathcal{B}}$ which converges to an element $W(Y)$ of $\af^l$. Using the notations in \ref{power_ident}, the sequence $(W_n)_{n\in\enn}$ converges to $W$ in $\aaa^l$.
Fix $d(X)=(d_{1}(X),\ldots,d_l(X))\in P$. We have
\begin{equation*}
  d(X)\circ W(Y)=\sum_{\alpha\in\enn^r}\sum_{\lambda=1}^ld_{\lambda\alpha}W_{\alpha}. \end{equation*}
But, the $d_{\lambda}(X)$'s being with finite support, there exists $N\in\enn$ such that $\alpha\leqslant_+\psi(N)$ for $\alpha$ in the union of the supports of the $d_{\lambda}(X)$'s. We then have
\begin{equation*}
  d(X)\circ W(Y)=\sum_{\alpha\leqslant\psi_r(N)}\sum_{\lambda=1}^ld_{\lambda\alpha}W_{\alpha}. \end{equation*}
Let  $V_N$ the element of the $0$-basis of $\aaa^l$ corresponding to $N$.  There exists $M\in\enn$ such that
\begin{align*}
  n\geqslant M&\Longrightarrow W_n-W\in V_N\\
  &\Longrightarrow\ W_{n\alpha}=W_\alpha\quad\text{for}\quad \alpha\leqslant_+\psi_r(N).
\end{align*}
It follows that
\begin{align*}
  d(X)\circ W(Y)&=\sum_{\alpha\leqslant_+\psi_r(N)}\sum_{\lambda=1}^ld_{\lambda\alpha}W_{n\alpha}\\ &=d(X)\circ W_n(Y)=0.
\end{align*}
Therefore, $W(Y)\in P^{\perp_{\af}}$ and $P^{\perp_{\af}}$ is closed in $\af^l$.\\

 $(3)\Longrightarrow (2)$: By \ref{orth-closed} $\boldsymbol{\mathcal{B}}=P^\perp$ for a subset $P$ of $\df^l$. But, $\boldsymbol{\mathcal{B}}$ being a submodule of $\af^l$, we know, by \ref{recipr-orths} that $P^\perp=P^{\perp_{\af}}$ So that $\boldsymbol{\mathcal{B}}=P^{\perp_{\af}}$.\qed
\begin{defi}[Systems] A  \textit{discrete algebraic dynamical system} of $\af^l$ (or simply \textit{system} of $\af^l$) is a
{subset of  $\af^l$} which verifies one of the equivalent conditions of Theorem \ref{autoregressive-orth-closed}.
\end{defi}
In order to characterize systems,  we introduce again more categorical notations: let $\Mod(\df)$ be the category of $\df$-modules. For two modules $M,N$ of $\Mod(\df)$, a morphism $f:M\longrightarrow N$ is an $\df$-homomorphism, thus an element of $\Hod(M,N)$. The subcategory of $\Mod(\df)$ whose objets are the finitely generated  is denoted by $\Modf(\df)$.
We will need the following lemma :
\begin{lem}[Finitely generated module, \cite{ob}]\label{finite-Mod}
Let be $M$ an objet of $\Modf(\df)$. Then, there exist integers $k,l\geqslant1$ and a matrix $R(X)\in\df^{k,l}$ such that the following isomorphism of $\df$-modules holds:
\begin{equation}
    \df^l/\df^k R(X)\cong M.\\
\end{equation}
\end{lem}

The element $M$ of $\Modf(\df)$ is then isomorphic to an element of $\Mod(\df)$ of the form $\df^l/P$ where $l\geqslant1$ and $P$ is a $\df$-submodule of  $\df^l$. So, we will have a closer look at these modules.\\

The following corollary is similar to \eqref{A-scal-prod}.
\begin{cor}[\cite{ob}]
For a $\df$-submodule $P\subset\df^l$, the map
\begin{align}\label{quot-scal-prod}
\begin{split}
\df^l/\langle P \rangle\times P^\perp &\longrightarrow\af\\
(\overline{d(X)},W(Y)) &\longmapsto d(X)\circ W(Y)
\end{split}
\end{align}
is again a scalar product. Thus we have the isomorphism

\begin{align}\label{quot-orth}
\begin{split}
P^\perp &\cong\Hod (\df^l/\langle P\rangle, \af)\\
W(Y) &\longmapsto
  \begin{cases}
   \langle -,W(Y)\rangle_{\af} : &\df^l/\langle P \rangle\longrightarrow \af \\
    &\overline{d(X)} \longmapsto d(X)\circ W(Y)
  \end{cases}
\end{split}
\end{align}
in $\Mod(\df)$.
\end{cor}

\begin{thm}[\cite{ob}]\label{system-hom}
A system $\boldsymbol{\mathcal{B}}$ is also on the form
\begin{equation}
    \Hod(M,\af)
\end{equation}
where $M\in \Modf(\df)$.
\end{thm}

\prf
Let $\boldsymbol{\mathcal{B}}=P^\perp$ be an system where $P\in\mathcal{M}_l.$ By \eqref{quot-orth}, we can write
\begin{equation*}
\boldsymbol{\mathcal{B}}=\Hod(\df^l/\langle P \rangle,\af)
\end{equation*}
with $\df^l/\langle P \rangle=M\in\Modf(\df)$. Conversely, if $M\in\Modf(\df)$, then, according to lemma  \ref{finite-Mod}, there exists $P\in\Modf(\df)$ such that we can write $M=\df^l/P$. Then, using again \eqref{quot-orth}, we have $\Hod(M,\af)=P^\perp$ and it is a system.\qed
\section*{Conclusion}
Main tools for our theory are the orthogonal $\perp$, constructed  with the scalar product $\langle -,-\rangle$ with value in  $\ef$ and the orthogonal $\perp_{\af}$, constructed with the scalar product $\langle -,-\rangle_{\af}$ with value in $\af$ and the functors $\Homf(-,\ef)$ and $\Hod(-,\af)$. We have constructed and used $0$-bases in the vector spaces $\aaa$  and $\aaa^l$ in our proof for  finding closed subspaces.
If $\subset\df^l$, then $P^\perp$ is a closed subspace of $\af^l$ and $P^{\perp_{\af}}$ is a closed submodule of $\af^l$.  The functors were respectively used to explain the polynomial-power series multiplication  ``$\circ$ '' and to characterize systems. 

\section*{Acknowledgements}
My special thanks  go to Ulrich Matthey, CEO of Kontecxt in Essen, Germany who has partially financed my research projects, with the ``Sur-Place'' Grant.

\end{document}